\documentstyle[psfig]{article}
\newcommand{\thsum}[1]{\;\#{\raisebox{-3pt}{$\ \!\!_{#1}$}}\,}
\newcommand{\ol}{\overline}

\newcommand{\R}{{\bf R}}

\newcommand{\C}{{\bf C}}

\newcommand{\p}{{\partial}}

\newcommand{\om}{{\omega}}
\newcommand{\eps}{{\varepsilon}}

\newcommand{\Ga}{{\Gamma}}

\newcommand{\io}{{\iota}}

\newcommand{\Ll}{{\cal L}}

\newcommand{\Nn}{{\cal N}}

\newcommand{\Tt}{{\bf T}}

\newcommand{\Tilde}{\widetilde}

\newcommand{\bcp}{{\Tilde{\C P}\,\!^2}}

\newcommand{\Si}{{\Sigma}}

\newcommand{\MS}{{\medskip}}

\newcommand{\NI}{{\noindent}}
\newcommand{\proof}[1]{\noindent{\bf Proof#1:\  }}

\newcommand{\QED}{\hfill$\Box$\medskip}

\newtheorem{theorem}{Theorem}[section]
\newtheorem{cor}[theorem]{Corollary}

\newtheorem{definition}[theorem]{Definition}
\newtheorem{example}[theorem]{Example}
\newtheorem{remark}[theorem]{Remark}
\newtheorem{lemma}[theorem]{Lemma}

\newtheorem{prop}[theorem]{Proposition}

\newcommand{\at}{{@}}

\title{Associativity properties of the symplectic sum}
\author{Dusa McDuff\thanks{Partially supported by
NSF grant DMS 9401443.} \\ State University of New York at Stony Brook \\
{\small ( dusa\at math.sunysb.edu)}\and Margaret Symington\\Stanford
University and Harvard University\\{\small (margaret\at math.stanford.edu)}}

\begin{document}
\maketitle
\begin{center} Preliminary version
\end{center}
\section*{Abstract}
In this note we apply a 4-fold sum operation  to
develop an associativity rule for the pairwise symplectic sum.  This
allows us to  show that certain diffeomorphic symplectic $4$-manifolds made out
of elliptic surfaces
are in fact symplectically deformation equivalent.    We also show that blow-up
points can be traded from one side of a symplectic sum to another without
changing the symplectic deformation class of the resulting manifold.

\section{Introduction}
Recently there have been several new constructions for compact symplectic
$4$-manifolds $(X, \om)$ as well as  great progress (via Taubes-Seiberg-Witten
theory) in understanding invariants for such manifolds.   One of the main
consequences of Taubes' work~\cite{TAU} is that the Gromov invariants of
$(X,\om)$ are invariants of the diffeomorphism type of $X$ rather than of  its
symplectomorphism type.
It would be very interesting to understand
whether or not a given diffeomorphism type can support two different
symplectic structures. In fact, as yet no $4$-dimensional
example is known of a compact manifold with two distinct structures,
though such examples were found by Ruan~\cite{RU} in dimensions $6$ and
higher.
The results presented here were developed to show that some possible
candidates  for such forms $\om, \om'$
are in fact   equivalent.

The appropriate notion of symplectic equivalence in the present context is that
of
weak deformation equivalence.  Specifically, two
symplectic forms $\om,\om'$ on $X$ are {\bf deformation equivalent}
 if there is a family of (possibly non-cohomologous)
symplectic forms $\om_t, 0\le t\le 1,$   such that   $ \om_0 =\om$ and $\om_1 =
\om'$, and two symplectic manifolds $(X,\om), (X',\om')$ are
{\bf weakly deformation
equivalent}  if there is a diffeomorphism $\phi:X\rightarrow X'$ such
that $\phi^*(\om') $  is deformation equivalent to $\om$.   For example, a
K\"ahler
manifold supports a well-defined deformation class of symplectic forms since
the
set of  K\"ahler forms compatible with a fixed complex structure is convex and
hence
path-connected.

Throughout this paper we restrict
to the $4$-dimensional case, though many of our results have higher dimensional
analogues.  When $S\subset X$ and $S'\subset X'$ are symplectically embedded
surfaces in the 4-manifolds $X,X'$, we write
$$
(X,S) = (X',S')
$$
 if there is a symplectomorphism from $X$ to $X'$ that takes
$S$ to $S'$, and
$$
(X,S) \cong (X',S')
$$
 if the manifold/submanifold pairs
are weakly deformation equivalent.  (This means that the forms $\phi^*(\om')$
and $\om$ are equivalent under a symplectic deformation $\om_t$ consisting
of symplectic forms which are nondegenerate on $S$.  For
example, this is always the case if $(X,\om_t)$ is K\"ahler and $S$ is a
complex
curve.)  Finally, by a {\bf triple} $(X,S,T)$ we mean a symplectic 4-manifold
$X$
with two symplectically embedded Riemann surfaces $S$ and $T$  which intersect
transverally with positive orientation in a single point.

In~\cite{G} Gompf developed a pairwise symplectic sum, observing that, when
 a pair of manifolds $X,X'$ are summed  along a pair of codimension two
submanifolds,  a transverse pair of submanifolds
can be summed at the same time provided that certain conditions are satisfied.
In
$4$-dimensions, the only pertinent condition  is  that the transverse
surfaces must have positive intersection with the submanifolds along which
the sum is being taken.    Indeed, consider triples
$$
(X_1, S_1, T_1),\quad   (X_2, S_2, T_2)
$$
for which
$$
g_{T_1} = g_{S_2},\quad \iota_{T_1} = -\iota_{S_2}, \quad \int_{T_1}\om_1 =
\int_{S_2}\om_2, $$
where $g_S$ denotes the genus of $S$ and $\iota_S$ is the Chern number of its
normal bundle, which in this setting is equal to
the self-intersection number of $S$.
Then one can form the pairwise sum of the manifold/submanifold
pairs $(X_1,S_1)$ and $(X_2,T_2)$ along the symplectomorphic
surfaces $T_1,S_2$:
$$
(X_1,S_1)\thsum{T_1=S_2}(X_2,T_2) = (X_1\thsum{T_1 = S_2}   X_2,\, S_1\#T_2)
$$
where $S_1\#T_2$ is the connected sum of surfaces isotopic to $S_1,T_2$.
This sum is described in detail in \S 2.
(Our notation in which  $T_1$ is glued to $S_2$
might seem a little awkward, but will prove to be very convenient.)

Our first observation is that the $4$-fold sum operation, which is
developed by the second author in~\cite{SYM}, is invariant under cyclic
permutations.  The $4$-fold sum is possible when four symplectic triples $$
(X_1, S_1, T_1),\quad   (X_2, S_2, T_2),\quad (X_3, S_3, T_3),\quad   (X_4,
S_4,
T_4)
$$
are such that each $(X_i, S_i, T_i)$ can be summed to
$  (X_{i+1}, S_{i+1}, T_{i+1}) $ along $T_i,S_{i+1}$
as above, where $i$ is understood mod 4.
To form the sum, remove all eight surfaces
$S_i, T_i$ and naively start making symplectic pairwise sums:
\begin{eqnarray*}
X_1\thsum{T_1 = S_2} X_2, & & \quad X_2\thsum{T_2 = S_3} X_3\\
X_3\thsum{T_3 = S_4} X_4& &\quad X_4\thsum{T_4 = S_1} X_1.
\end{eqnarray*}
We explain in \S2 how these sums continue into the neighborhoods of the
intersection points to yield a smooth symplectic manifold, the $4$-fold sum.
In \S~3 we prove:

\begin{prop}[4-fold sum rule]\label{prop:4sum}
If triples $(X_i, S_i, T_i), 1\le i\le 4$ are such that for all $i\  (mod\,4)$
$$
g_{T_i} = g_{S_{i+1}},\quad
\iota_{T_i} = -\iota_{S_{i+1}},\quad
\int_{T_i}\om_i = \int_{S_{i+1}}\om_{i+1} $$
then
\begin{eqnarray*}
& & ( X_1 \thsum{T_1 = S_2} X_2 )
\thsum{S_1\#T_2=S_3\#T_4}
( X_3 \thsum{T_3 = S_4} X_4 ) \\
& & \qquad \qquad \qquad=
( X_4 \thsum{T_4 = S_1} X_1 )
\thsum{S_4\#T_1=S_2\#T_3}
( X_2 \thsum{T_2 = S_3} X_3 ).
\end{eqnarray*}
\end{prop}

Using this we prove an associativity rule for a sum of three triples
$(X_i, S_i, T_i)$, $i = 1,2,3$.  Before stating this, we
note that if $S$ and $T$ are symplectic submanifolds
which intersect positively along a symplectic submanifold
it is always possible to construct a symplectic submanifold
in the class $[S]+[T]$ which equals $S\cup T$ except near the intersection
$S\cap T$.  We will
think of this manifold as the  desingularization of $S\cup T$ and will denote
it by
$S+T$.
It is not hard to show that $S+T$ is well-defined up
to symplectic isotopy.

\begin{prop}[Associativity]\label{prop:assoc}  Suppose that for
$i=1,2,3\ (mod\, 3)$
\begin{eqnarray*}
g_{T_i} = g_{S_{i+1}},& & \int_{T_i}\om_i = \int_{S_{i+1}}\om_{i+1},\qquad
\mbox{ and
 }\\
\io_{T_1} = -\io_{S_{2}}, \quad \io_{T_2} & = & -\io_{S_{3}},\quad
 \io_{T_3} = -(\io_{S_1} + 2).
\end{eqnarray*}
Then
$$
(X_1\thsum{T_1 = S_2}X_2)\thsum{S_1\#T_2 = S_3 + T_3} X_3\;\; =\;\;
X_1\thsum{S_1+T_1=S_2\#T_3}(X_2\thsum{T_2=S_3}X_3).
$$
\end{prop}
\proof{} (Sketch)  To prove the result up to deformation we will find a
triple $(X_4, S_4, T_4)$ which, when summed with either
$X_3$ or $X_1$, yields a manifold symplectically deformation equivalent to
the original summand
and containing $S_3\# T_4$ or $S_4\#T_1$ as a representative
of the class $[S_3]+[T_3]$ or $[S_1]+[T_1]$ respectively.

The triple we will use is $(W_{g}, \Ga_{-k+2}, \Ga_k)$ where $W_g$ is
a ruled surface over a Riemann surface of genus~$g$ with a pair of sections
$\Ga_{-k+2},\Ga_k$ that intersect positively in one point.  Indeed, we
prove in \S~3:

\begin{lemma}\label{le:resol}
Given any triple $(X,S,T)$ and sufficiently small $\eps$, if
$(g,k)=(g_S,-\io_S)$ then the symplectic sum with the triple
$(W_g,\Ga_{-k+2},\Ga_k)$ yields
$$
( X , S+T) \cong ( W_{g}\thsum{ \Ga_k=S}X ,\, \Ga_{-k+2}\# T).
$$
Analogously, if $(g,k) = (g_T , \io_T +2)$ then
$$
( X , S+T) \cong ( X\thsum{T = \Ga_{-k+2}} W_{g},\, S\#\Ga_k).
$$
\end{lemma}

Granted this, we  apply the $4$-fold sum rule with
 $(X_4, S_4, T_4)= (W_g,\Ga_{-k+2},\Ga_k)$
where
$$
g = g_{S_1} = g_{T_3},\quad k = -\io_{S_1} = \io_{T_3}+2.
$$
It is easy to check that the $4$-fold sum is  well-defined.  Moreover,
by Lemma~\ref{le:resol}
\begin{eqnarray*}
(X_4\thsum {T_4=S_1 } X_1,\, S_4\#T_1) & \cong & (X_1, S_1 + T_1)
\quad\mbox{and}\\
 ( X_3 \thsum{T_3=S_4} X_4,\, S_3\#T_4) & \cong &  ( X_3, S_3+T_3).
\end{eqnarray*}
This proves the two manifolds are deformation equivalent.  The proof that
they are symplectomorphic is given in \S4. \QED

\begin{remark} \rm
In order to prove Proposition~\ref{prop:assoc} up to
symplectomorphism rather than deformation equivalence we will need a
thickening/thinning procedure which is described in \S4.
By thickening and thinning we can build the neccessary $X_4=W_g$ out
of pieces removed from the other $X_i$'s.
\end{remark}

\MS

One application of these results is to prove the symplectic equivalence of two
manifolds which are constructed out of elliptic surfaces.  Let $E(n)$ be the
elliptic
surface which is the $n$-fold branched cover of
$E(1) = \C P^2\# 9\ol{\C P}^2$ along a
fiber.  Then $E(n)$ contains  $9$ sections $\Si_{-n}$ which are $n$-fold covers
of
the exceptional spheres in  $E(1)$ and so have self-intersection $-n$.  In
terms of
the symplectic sum, we can inductively define
$$
(E(n), \Si_{-n}) \cong (E(n-1)\thsum{F_{n-1} = F_1} E(1),\, \Si_{-n+1}\#
\Si_{-1})
$$
where $F_k$ is a generic fiber in $E(k)$.
When $n = 4$ one can sum along
the section $\Si_{-4}$ and a quadric curve $Q$ in $\C P^2$ to form the
manifold $E(4)\thsum{\Si_{-4} = Q}  \C P^2$ which is not diffeomorphic to any
complex surface (see Gompf [1]).
On the other hand,
when $n = 3$ there is a torus $T_{-1}$ in the homology class
$[\Si_{-3}] + [F_3]$.
Since
$$
T_{-1}\cdot T_{-1} = \Si_{-3}\cdot \Si_{-3} + 2 \Si_{-3}\cdot F_3 = -1,
$$
 $E(3)$
can be summed with $Y =  \C P^2\# 8 \ol{\C P}^2$ along the tori
$T_{-1},T_1$ where $T_1$ is a torus of self-intersection $+1$ in $Y$ which is
obtained from an elliptic curve in $\C P^2$ through the $8$ blown up points.

Stipsicz~\cite{STI} proved that the Donaldson invariants of these
two manifolds
$$
E(4)\thsum{\Si_{-4} = Q}  \C P^2, \quad E(3)\thsum{T_{-1} = T_1} Y
$$
are the same and Gompf~\cite{G2} has shown
using Kirby calculus that they are diffeomorphic.
We show here that the manifolds
are symplectically deformation equivalent.
(Observe that, since they are built out of K\"ahler surfaces,
both manifolds have symplectic forms well-defined up to deformation.)

\begin{prop}\label{prop:main} The manifolds
$ E(4)\!\thsum{\Si_{-4} = Q}\!\C P^2$
and $ E(3)\!\thsum{T_{-1} = T_1}\!Y$ are symplectically deformation
equivalent.   \end{prop}
\proof{}
First, scale the symplectic forms on $E(1),E(3)$ so that the symplectic areas
of the
fibers are equal.  Next, adjust these forms
(by pulling back suitable forms from
the base of the elliptic fibrations) to make the sections
$\Si_{-1}, \Si_{-3}$ have
the same symplectic area $a_\Si$, and choose a symplectic form on $\C P^2$ such
that the symplectic area of a line equals $a_\Si$    Then,  if we take
\begin{eqnarray*}
(X_1, S_1,T_1) & = & (E(3), \Si_{-3}, F_3)\\
(X_2, S_2,T_2)  & = & (E(1),  F_1,\Si_{-1})\\
(X_3, S_3,T_3) & = & (\C P^2, L_1, L_2)
\end{eqnarray*}
where $L_1, L_2$ are two lines in $\C P^2$, the $3$-fold sums of
 Proposition~\ref{prop:assoc} are defined.  Further,
\begin{eqnarray*}
(X_1\thsum{T_1 = S_2}X_2)\thsum{S_1\#T_2 = S_3 + T_3} X_3
& = & (E(3)\thsum{F_3 = F_1} E(1))\thsum{\Si_{-4} = Q}\,\C P^2\\
& = & E(4) \thsum{\Si_{-4} =Q}\,\C P^2.
\end{eqnarray*}
On the other hand
\begin{eqnarray*}
X_1\thsum{S_1+T_1=S_2\#T_3}(X_2\thsum{T_2=S_3}X_3)
 & = &
E(3)\thsum{T_{-1} = F_1\#L_2} (E(1)\thsum{\Si_{-1} = L_1 } \C P^2)\\
 & = &  E(3)\thsum{T_{-1} = T_1} Y,
\end{eqnarray*}
where the last equivalence holds because the  sum with $\C P^2$ is
just a symplectic blow-down of $\Si_{-1}$ and so takes $E(1)$ to $Y$
and $F_1$ to $T_1$. \QED

Another application in the same spirit uses properties of the ruled surfaces
$W_{g}$ to show that blow-up points can be traded
from one side of a symplectic sum to
the other without changing the deformation class of the symplectic structure.
McCarthy and Wolfson~\cite{MW} noted that a standard handle trading argument
shows that this can be done up to diffeomorphism, as explained in detail
by Gompf in Lemma 5.1 of~\cite{G}.  However, these authors left
open the question of symplectic equivalence.  We prove

\begin{prop}\label{prop:bl}
Consider  symplectic  pairs $(X,S), (X',S')$  such that
$$
\io_S = - \io_{S'} + 1,\quad g_S = g_{S'}.
$$
Let $\Tilde X$ be the blow-up of $X$ at a point of $S$ and $\Tilde S$ the
proper
transform of $S$, and similarly for $(\Tilde { X'}, \Tilde {S'})$.  Then
 $$
X \,\# _{S = \Tilde {S'}}\,\Tilde {X'} \;\;\cong \;\;\Tilde X
\,\# _{\Tilde S=S'}\,X'.
$$
\end{prop}

\begin{remark}\rm (i)
Because the blow-down operation on $\Tilde X$ may be interpreted as a sum
with $\C P^2$, an equivalent way of stating this (which shows its similarity to
the
associativity rule) is:
$$
(\Tilde X \thsum{E = L_1} \C P^2)
\,\# _{\Tilde{S}\#L_2 = \Tilde {S'}}\,\Tilde {X'}
\;\;\cong\;\;  \Tilde {X} \,\# _{\Tilde {S} = L_2\#\Tilde{S'}}
\,(\C P^2 \thsum{L_1 = E'} \Tilde {X'}).
$$
Here $E, E'$ denote the exceptional spheres in $\Tilde X, \Tilde {X'}$, and
$L_1,
L_2$ are lines in $\C P^2$.

\MS

\NI
(ii) The invariance of the symplectic structure under the trading
of blow-up points can be at
most up to deformation equivalence since it is impossible to fix the symplectic
areas
of $S$ and $S'$ in such a way that both sums can be performed.
To see this, observe that the area of the proper transform of a surface (after
a
symplectic blow-up) is less than that of the original
surface, so the symplectic sum along $S,\Tilde{S'}$ requires
$\int_S\om < \int_{S'}\om'$ while a sum along $\Tilde{S},S'$ requires the
reverse inequality.
\end{remark}

As another
application, we show that rational blow-down of a $-4$-sphere gives nothing new
if this sphere is the blow-up of a $-3$-sphere. Again, this
was proved by Gompf as far as concerns diffeomorphism type.

\begin{cor}\label{Q}  Let $S\subset X$  be a sphere with $\io_S = -3$ and let
$Q$
be a quadric curve in $\C P^2$.  Then
$$
\Tilde X\,\# _{\Tilde S = Q}\,\C P^2 \cong X.
$$
\end{cor}
\proof{} Proposition~\ref{prop:bl} shows that
$\Tilde X\# _{\Tilde S = Q}\C P^2
\cong X\# _{ S = \Tilde Q}\,\bcp$.  But, we may think of  $\bcp$ as the
projectivization ${\bf P}(\C \oplus \Ll_3)$ of a complex
rank $2$ bundle over $\C P^1$
where $\Ll_3$ is a complex line bundle of Chern number $3$.  Thus $\bcp$  is
the
union of a neighborhood of the section $\Tilde Q$ with a neighborhood of a $-3$
section, and it follows immediately that $X\# _{ S = \Tilde Q}\,\bcp \cong X$.
\QED

\NI
{\it Acknowledgement}
The first author wishes to thank Gompf for bringing these examples to her
attention and the organizers of the G\"okova conference for providing such a
beautiful place in which to do research.  The second author wishes to
thank Yakov Eliashberg for his guidance in thinking about symplectic sums
along intersecting submanifolds.

\section{The $4$-fold  sum}

The associativity rule (Proposition~\ref{prop:assoc})
is a consequence of the fact that a simple version
of the 4-fold sum is equivalent to a sequence of three symplectic
sums (two of which are pairwise) performed in either of two ways.
We begin with a description of the pairwise symplectic
sum in terms of images under the moment map for a local torus action.

\subsection{The Symplectic and Pairwise Sums}

Given a pair of symplectic submanifolds $S_i \subset (X_i,\om_i)$, $i=1,2$
and a symplectomorphism $\phi:S_1\rightarrow S_2$, one can  perform a
symplectic sum of $X_1$ and $X_2$ along $S_1$ and $S_2$ provided
the normal numbers (Chern numbers of the normal bundles) of the submanifolds
sum to zero.  A good way to see this operation is as an inverse to Lerman's
symplectic cutting
procedure~\cite{LER}.  Observe first that a
codimension two symplectic submanifold
$S$ in $X$ always has a tubular neighborhood $\Nn_S$
that admits a Hamiltonian circle action with
fixed point set  $S$.   Moreover, one can clearly extend the induced action on
$\Nn_S - S$ to a free Hamiltonian action on a collar neighborhood
 $\overline\Nn_S$
of the boundary in an appropriate closure
$\overline{X-S}$ of $X - S$.  Here  $\Nn_S$ is an open disc bundle over $S$,
and $\overline \Nn_S$ is the associated bundle with fibers $[0,1)\times S^1$
so that its boundary $\p_S = \{0\}\times S^1$ is a circle bundle over $S$.
 Furthermore, this
boundary $\p_S$ is a level set of the Hamiltonian which generates the
free action, and its symplectic
reduction is $S$ itself.  Now, the way to form the
sum
$$
X_1 \thsum{S_1=S_2} X_2
$$
is to remove the
submanifolds $S_i$, take the closures
$\overline{X_i - S_i}$ as described above,
and then identify the boundaries  via an orientation
reversing diffeomorphism $\overline\phi: \p_{S_1}\to \p_{S_2}$ that covers
$\phi$ and thus matches the  characteristic foliations
(along which the symplectic
forms are degenerate).

\begin{remark}\rm
Note that when the normal bundles of the submanifolds are trivial, the
diffeomorphism class of the summed manifold depends on the choice
of the fiber isotopy class of the map $\overline\phi$.  In the examples
we consider in this paper we sum along fibers in elliptic surfaces
and use the canonical framings to get the boundary identifications.
\end{remark}

In order to describe the pairwise sum in a similar way, one needs to use
torus actions rather than circle actions.  Recall that when a closed
4-manifold admits a Hamiltonian  action of $T^2$,
the image of an associated moment map is a convex polytope in ${\R}^2$.
The preimages of points on the interior of the polytope are tori,
while the preimages of points on an edge or vertices are circles or
points respectively.  We use the convention that replacing a solid
line segment
in the image of a moment map by a heavy dotted
line segment corresponds to replacing the submanifold $S$
(the preimage of the solid line)
with the  associated boundary $\p_S$.  A light dotted line
indicates an
open boundary.  For instance,

\begin{figure}
\centerline{\psfig{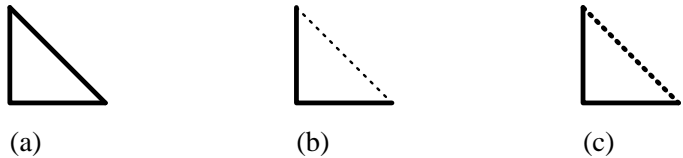}}
\caption{Images of moment maps.\label{fig:cp2}}
\end{figure}

\begin{example}\rm
Figure~1 shows the
image under the moment map for (a) $\C P^2$, (b) the open
ball obtained by removing a $\C P^1 = S$ and
(c) its closure, a closed ball with boundary $\p_S$.
\end{example}

For the purposes of this paper we want to keep track of the effect of
the symplectic sum on transverse symplectic submanifolds.
If two symplectic surfaces $S,T$ intersect transversely and positively, then
one of them can be perturbed, via an isotopy through
symplectic surfaces, so that the intersection is orthogonal with respect to the
symplectic structure (see~\cite{G} Lemma 2.3, for example).

Consider a triple $(X,S,T)$ as in \S 1.
Then in a neighborhood $\Nn_x$ of the intersection point $\{x\} = S \cap T$
there is a Hamiltonian
$T^2$ action such that the first $S^1$ factor has fixed point set
$S\cap\Nn_x$ and the
second has fixed point set $T\cap\Nn_x$.  Thus the image of the moment map is a
neighborhood of a corner in a square.
We may close $X - (S\cup T)$ by
adding a boundary (with corner) to get a compact symplectic manifold on which
there is a
free local $T^2$
action acting in a neighborhood of the corner.

For triples $(X_i,S_i,T_i)$, $i=1\cdots 4$ we
choose neighborhoods of the intersection points
whose images under the moment map are as
in Figure~2(a), where the slopes of the
slanted edges are $-\io_{S_i}$ and $-1/\io_{T_i}$.

\begin{remark}\rm
These images show the correct
convexity (or concavity) of the symplectic neighborhoods.
Indeed, when the normal number of a submanifold is positive,
removing a tubular neighborhood whose boundary is a level set of the
Hamiltonian and taking the symplectic reduction of the newly formed
boundary yields a surface whose area is smaller than that of the original
surface.
\end{remark}

In order for our notation to be consistent with that used
for the $4$-fold sum, the summing operation identifies
$T_i$ with $S_{i+1}$ (mod 4).  We assume that $S_i$ intersects
$T_i$ orthogonally in $x_i$ and that the gluing map
$
\phi_i: { T_i}\to {S_{i+1}}
$
takes $x_i$ to $ x_{i+1}$.
Then the first diagram of
Figure~2(b) shows the image under the moment map
of these neighborhoods after the  symplectic sum has been taken along
$T_1,S_2$  using the symplectomorphism $\phi_1$.
The bold horizontal line at the bottom  consists of points
with preimage equal to a circle and so  is a neighborhood
in the connected sum $S_1\# T_2$ of the attaching circle.
The heavy dotted vertical line segment represents the normal
$2$-disc bundle over the attaching circle, or equivalently,
the intersection of the neighborhood $\Nn_x$ and the identified boundaries
associated to the submanifolds $T_1,S_2$.  The second diagram
is a similar picture of the sum along $T_3, S_4$.
Since any positive intersection can be made orthogonal via an isotopy
of one of the intersecting surfaces,
Figure 2(b) makes the following lemma clear.
It is a rephrasing of the $4$-dimensional  case of Theorem~1.4 in
Gompf~\cite{G}.

\begin{lemma}\label{lemma:sum}
Consider two triples $(X_i,S_i, T_i)$, $i=1,2$.
If $T_1$ and $S_2$ have the same area and genus, and $\io_{T_1}=-\io_{S_2}$,
then in the manifold
$X_1
\thsum{T_1=S_2} X_2$ there is a symplectic surface that is
the connected sum of surfaces isotopic to $S_1$ and $T_2$.  The normal
number of the surface $S_1 \# T_2$
is the sum of the normal numbers of $S_1$ and $T_2$.
\end{lemma}

\subsection{The 4-fold Sum}

The 4-fold sum is a sum along the four pairs of surfaces in
four triples $(X_i,S_i,T_i)$, $i=1,\cdots 4$
such that for each $i$,
\begin{itemize}
\item $S_i\cap T_i= \{x_i\}$ and the intersection is orthogonal with respect
to $\omega_i$,
\item $\int_{T_i} \omega_i = \int_{S_{i+1}} \omega_{i+1}$ and
\item $\io_{T_i} = - \io_{S_{i+1}} $
\end{itemize}
where the subscripts are taken mod 4.  Call such a collection of triples
{\bf admissible}.

Because the normal numbers of each pair of symplectomorphic surfaces sum
to zero, it is possible  to
sum along all four pairs $T_i,S_{i+1}$.  Trying to do these sums
simultaneously leads to the following definition:

\begin{definition}\rm
Given an admissible collection of triples $(X_i,S_i,T_i)$,
$i=1,\cdots 4$ and symplectomorphisms $\phi_i:T_i\rightarrow S_{i+1}$
such that $\phi_i(x_i)=x_{i+1}$ (mod 4),
let
$$
\overline{X_i-(S_i\cup T_i)}
$$
 be the closure of $X_i-(S_i\cup T_i)$
with free local
$T^2$ action as described above.  Then, choosing
orientation reversing diffeomorphisms $\ol\phi_i$ that cover the $\phi_i$,
we define
$$
\sqcup_{i=1}^4\;\;
 \overline{X_i-(S_i\cup T_i)}\;/\;\overline\phi_i
$$
to be the {\bf 4-fold sum} of the $X_i$ along the surfaces $T_i,S_{i+1}$.
\end{definition}

\begin{figure}
\centerline{\psfig{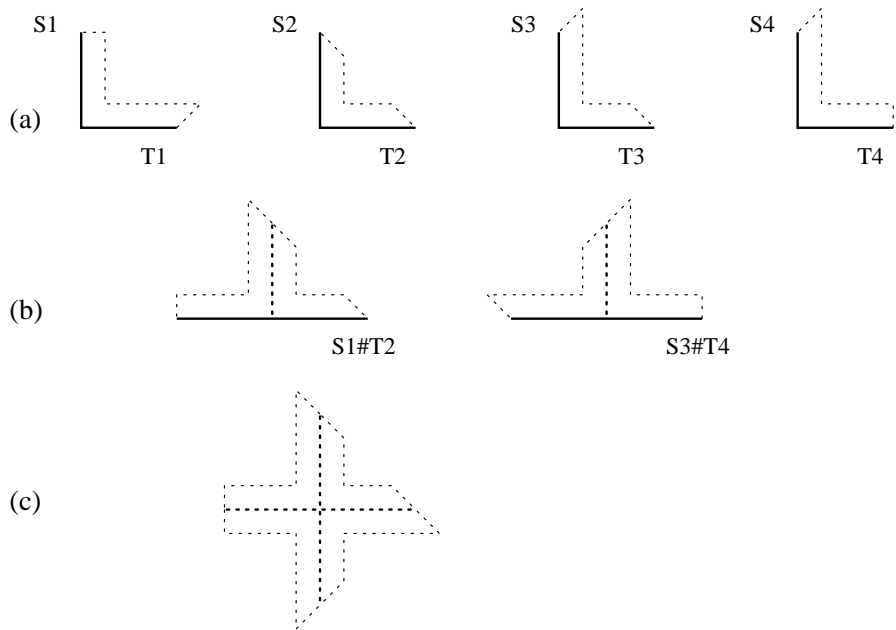}}
\caption{The 4-fold sum.\label{fig:4fold}}
\end{figure}

To see that this construction yields a smooth symplectic manifold,
notice that it is equivalent to a sequence of three symplectic  sums.
Indeed, given a set of four triples as in the definition,
by Lemma \ref{lemma:sum} we can use the maps $\overline\phi_1$ and
$\overline\phi_3$ to form
pairwise sums along the pairs $T_1,S_2$ and $T_3,S_4$
to yield two manifolds which contain surfaces $S_1 \# T_2$ and
$S_3 \# T_4$ respectively.
These surfaces have the same area and genus, and have normal numbers
$\io_{S_1 }+ \io_{T_2}$ and $\io_{S_3} + \io_{T_4} = -\io_{T_2} - \io_{S_1}$
respectively.
Therefore we can form the symplectic sum along this new pair.  In fact, to
perform this sum, we can use the diffeomorphisms $\overline\phi_2$ and
$\overline\phi_4$ which agree on the overlapping circle.
Figure~2 shows these sums, keeping track of the convexity of
all of the tubular neighborhoods.
Thus, the only place the 4-fold sum might not have been a smooth symplectic
manifold it is in fact symplectomorphic to the product of $T^2$ and a domain in
${\bf R}^2$ with the symplectic structure of $T^*T^2$.
A precise statement and proof of this fact is given in Symington~\cite{SYM}.

Proposition~\ref{prop:4sum} follows immediately since the 4-fold sum is
certainly also equivalent to first summing pairwise along
the surfaces $T_2,S_3$ and $T_4,S_1$ and then summing the resulting
manifolds along the surfaces $S_2\#T_3, S_4\#T_1$.

\section{Proofs of the main results}

We begin by describing the manifolds $W_g$ of Lemma~\ref{le:resol}.
For each genus $g$ and integer $n\ge 1$
let $W_g = W_{g,n}$ be an $S^2$ bundle over a Riemann surface of
genus~$g$.  Make $W_g$ be the trivial bundle (i.e. a product) if $n$ is
even, and the nontrivial bundle when $n$ is odd.
In either case, there
is a unique symplectic
structure on the ruled surface $W_g$ up to deformation.
(In fact, by~\cite{LM} symplectic forms
in a given cohomology class on $W_g$ are unique up to isotopy.)
We choose a symplectic structure on $W_g$ so that the manifold
contains symplectic sections $\Ga_{-n+2p}$ of self-intersection $-n+2p$
for all $p \ge 0$.
Then in particular there are symplectic sections $\Ga_{-k+2}$ and
$\Ga_{k}$ which intersect once transversally (and positively) for each
$-n\le k\le n+2$.

To be explicit, realize $W_g$ as
the projectivization of the complex rank $2$ bundle $\C \oplus \Ll_{n}$
as in Corollary~\ref{Q}.  Then $\Ga_n$ and
$\Ga_{-n}$
are the holomorphic (and therefore symplectic) sections at zero and infinity.
Observe that they are disjoint.
If $F$ is a fiber then $[\Ga_{-n}]$ and $[F]$ form
a basis for the homology of $W_g$ such that $[\Ga_n] = [\Ga_{-n}] + n[F]$.
Because the fiber can also be realized as a holomorphic curve it is clear
that there are also holomorphic (and hence symplectic) sections
$\Ga_{-n+2} \in [\Ga_{-n}]+[F]$ and  $\Ga_{n-2} \in [\Ga_{-n}]+(n-1)[F]$.
Notice that the ruled surface
$W_g$ is also the total space of a Hamiltonian
$S^1$-action with fixed point sets $\Ga_n$ and $\Ga_{-n}$, or indeed of
a Hamiltonian with fixed point sets $\Ga_{n-p}$ and $\Ga_{-n+p}$ for
any $0\le p\le n-1$.

Henceforth, we assume that $W_g=W_{g,n}$ with $n$ sufficiently large
that any sections referred to can be (and are) taken to be symplectic.
Note that the parity of the self-intersection numbers of the sections
will in all cases make it clear whether or not $W_g$ is the trivial bundle.

We claim that
the manifold $W_g$ has the exact properties we need in order to prove
Lemma~\ref{le:resol}.  First of all,
because the sections of opposite self-intersection are disjoint,
$[\Ga_k] \cdot ([\Ga_{-k}]+[F]) = 1$, so
the sections $\Ga_k$ and $\Ga_{-k+2}$ intersect
once positively.
Second of all, we can view $W_g$
either as the union of symplectic neighborhoods
of $\Ga_k$ and $\Ga_{-k}$ or of $\Ga_{k-2}$ and $\Ga_{-k+2}$, thanks to
the fibration by symplectic spheres.

We now show that summing $(W_{g}, \Ga_{-k+2}, \Ga_k)$ with a triple
$(X,S,T)$
desingularizes the intersecting submanifolds $S,T$.
\MS

\NI
{\bf Proof of Lemma~\ref{le:resol}}  We will show that if
 $S\subset X$ has genus $g$ and $\io_S = -k$,  then
$$
( X, S+T) \cong  ( W_{g}\thsum{ \Ga_k=S}X ,\,\Ga_{-k+2}\# T).
$$
The other identity then follows by replacing $k$ by $-k+2$ and interchanging
the roles of $S,T$.

Now, it is obvious that $ X \cong W_{g}\!\thsum{ \Ga_k=S}X
$, i.e. $W_{g}$ is a neutral element for the symplectic sum operation
in the category of symplectic deformation equivalence classes.
As for the statement about the submanifold $S+T$, observe that
we have
$$
(X, S, T) \cong ( W_g\thsum{ \Ga_k=S}X,\, \Ga_{-k},\, F\# T).
$$
Hence
\begin{eqnarray*}
(X, S+T) & \cong &  (W_g\thsum{ \Ga_k=S}X, \,\,\Ga_{-k} + F\#T) \\
& \cong &  (W_g\thsum{ \Ga_k=S}X,\,\, (\Ga_{-k}+F)\# T )\\
& \cong &   (W_g\thsum{ \Ga_k=S}X, \,\,\Ga_{-k+2}\# T),
\end{eqnarray*}
where the second equivalence  holds because $\Ga_{-k}$ is disjoint from $
\Ga_{k}$. \QED

To show that one can trade blow-up points we use a pair of disjoint
sections in the ruled surface $W_g$, rather than a pair that intersect
once.
\MS

\NI
{\bf {Proof of Proposition~\ref{prop:bl}}}
Let $W = W_{g}$ be a ruled surface with disjoint sections $\Ga_{\pm(k+1)}$
and let $W' = W'_{g}$ be one with disjoint sections $\Ga_{\pm k}$.  Thus, one
of these manifolds $W,W'$ will be a trivial fibration and the other one
non-trivial.
Blow $W$ up at a point of
$\Ga_{k+1}$ and blow $W'$ up at a point of $\Ga_{-k}'$. It is well
known that there is a diffeomorphism that realizes $\Tilde W \cong\Tilde{W'}$,
taking $\Tilde{\Ga}_{k+1} $ to $\Ga_k'$, and $\Ga_{-k-1}$
to $\Tilde{\Ga'}_{-k}$.  Recall that the symplectic structures on these
manifolds are equivalent up to deformation.
Thus
$$
( \Tilde W, \Tilde{\Ga}_{k+1}, \Ga_{-k-1}) \cong (\Tilde {W'}, \Ga_k',
\Tilde{\Ga'}_{-k}).
$$

Letting $k=-\io_S-1$, we can deform
$(X,S) \cong (X\!\thsum{S = \Ga_{-k-1}}\!W, \Ga_{k+1})$
and we clearly have:
\begin{eqnarray*}
(\Tilde X, \Tilde S)  & \cong&  (X\thsum{S = \Ga_{-k-1}} \Tilde W,
\Tilde{\Ga}_{k+1})\\ & \cong &
(X\,\# _{S= \Tilde{\Ga'}_{-k}}\, \Tilde{W'}, \Ga_k').
\end{eqnarray*}
Hence, because $\Tilde{\Ga'}_{-k}$ and $\Ga_k'$ are disjoint in $\Tilde{W'}$,
\begin{eqnarray*}
\Tilde X\,\# _{ \Tilde S = S'}\, X' & \cong &
(X\,\# _{S = \Tilde{\Ga'}_{-k}}\, \Tilde{W'})\thsum{\Ga_k' = S'} X'\\
& \cong & X\,\# _{S = \Tilde{\Ga'}_{-k}}\,(\Tilde{W'}\thsum{\Ga_k' = S'} X')\\
& \cong & X\,\# _{S = \Tilde{S'}}\, \Tilde X'.
\end{eqnarray*}
Note that the effect of deforming $(X,S) $ to
$ (X\!\thsum{S = \Ga_{-k-1}}\!W, \,\Ga_{k+1})$ is to
localize the argument near $S$, i.e. we represent a neighborhood of $S$ in
$X$ as a
neighborhood of $\Ga_{k+1}$ in $W$ and then  work in $W$.   \QED

\section{Thickening and thinning}

In order to prove the
associativity rule up to symplectomorphism we  need a  refinement of
Lemma~\ref{le:resol}.   The notation $W_{g,\eps}$ will mean that the
symplectic form on $W_{g}$ has been scaled so that the fiber has symplectic
area $\eps>
0$.  By~\cite{LM} the symplectic form on $W_{g,\eps}$ is then determined up to
symplectomorphism (even isotopy) by specifying the area of one  section.
(The only condition on this area is the following: if $\om(\Ga_k) = a$
we need $a > k\eps/2$, unless $g = 0$, $k$ is odd, in which case we need $a >
(k+1)\eps/2$. Thus, given $a$, this is satisfied for sufficiently
small $\eps$.)

 The next lemma says that $W_{2\eps}$ can be thought of as the sum of two
copies of
$W_\eps$.

\begin{lemma}\label{le:eps0} If  symplectic forms on two copies  $W_{g,\eps}^i,
i = 1,2$
of $W_{g,\eps}$ are chosen  so that the area of the section
$\Ga_k^1$ in the first equals the
area of $\Ga_{-k}^2$ in the second, then
$$
(W_{g,\eps}^1\thsum{\Ga_k^1 = \Ga_{-k}^2}W_{g,\eps}^2, \Ga_{-k}^{1}, \Ga_{k}^2)
=
(W_{g,2\eps}, \Ga_{-k},\Ga_{k}).
$$
 \end{lemma}
\proof{}  We use Lerman's symplectic cutting procedure~\cite{LER}
to show that $W_{g,2\eps}$
decomposes as a sum of this kind. Choose a Hamiltonian function $f$ on
$W_{g,2\eps}$
with fixed point sets $\Ga_{-k}=f^{-1}(0)$ and $\Ga_{k}=f^{-1}(2\eps)$.
Then both of the ruled manifolds obtained by
cutting $W_{g,2\eps}$ along the $S^1$-invariant hypersuface $f^{-1}(\eps)$
and taking the $S^1$ reduction
along the boundaries of $f^{-1}([0,\eps]),f^{-1}([\eps,2\eps])$
have fiber of size $\eps$ and so may be identified with the manifolds
$W_{g,\eps}^1,W_{g,\eps}^2$.
\QED

Given a triple $(X,S,T)$ with symplectic structure $\om$
we write
$$
(\Tt_S^{-}(X), S^-, T^-)
$$
for the
(deformation equivalent)  manifold formed by \lq\lq thinning\rq\rq\,
$X$ along $S$ by the amount $\eps$, for some sufficiently small $\eps$
(which determines $\Tt_S^{-}(X)$ up to symplectomorphism).
In terms of the language of \S2,
we remove an $S^1$-invariant open
tubular neighborhood $\Nn_S$ of $S$ with fiber of area $\eps$,
and then reduce the boundary of $X-\Nn_S$ by the $S^1$ action.
The surface $S^-$ is the
symplectic reduction of the boundary.
Because $T$ coincides near $S$ with an orthogonal symplectic fiber of $\Nn_S$,
the manifold
$T^-$ is just $T$ with a disk of area $\eps$ removed and the circle boundary
collapsed to a point.  Indeed, we have
$$
(X,S,T) = (\Tt_S^{-}(X) \thsum{S^- = \Ga_k} W_{g,\eps},\, \Ga_{-k},\, T^-\#F)
$$
where $g=g_S$, $k=-\io_S$ and $F$ is a fiber of the ruled surface $W_{g,\eps}$.
Letting $\om^-$ be the symplectic form on $\Tt_S^-(X)$, we have
$$
\int_{\Tt_S^{-}(X)} (\om^-)^2 < \int_X\om^2,\qquad
\int_{S^-}\om^- = -\eps\io_S + \int_S\om, \qquad
\int_{T^-}\om^- = -\eps + \int_T\om.
$$

Analogously, we can \lq\lq thicken\rq\rq\,
along the surface $S$ by removing $S$, taking the closure $\ol{X-S}$ as
in \S~2, and
gluing in an $S^1$ invariant neighborhood of a surface $S^+$
diffeomorphic to $S$.  Again take the area of the fibers to be
some sufficiently small $\eps$.
We denote the triple that arises from the thickening by $\eps$ along
$S$ by
$$
(\Tt_S^{+}(X),S^+,T^+).
$$
It is easy to see that this thickening is just given by summing with
$W_{g,\eps}$:
$$
(\Tt_S^{+}(X),S^+,T^+) = (X \thsum{S = \Ga_k} W_{g,\eps},\, \Ga_{-k},\, F\#T).
$$
Moreover, the symplectic structure $\om^+$ on this manifold has
$$
\int_{\Tt_S^{+}(X)} (\om^+)^2 > \int_X\om^2,\qquad
\int_{S^+}\om^+ = \eps\io_S + \int_S\om, \qquad
\int_{T^+}\om^+ = \eps + \int_T\om.
$$

We will use the following property of thickening and thinning in the proof
of the associativity rule.

\begin{lemma}\label{le:eps}
If $(X_i,S_i,T_i), i = 1,2$ are triples such that
the sum
$$
\left(X_1\thsum{T_1 = S_2} X_2, S_1\# T_2\right)
$$
is defined, then for sufficiently small $\eps$
$$
\left(\Tt_{T_1}^-(X_1) \thsum{T_1^- = S_2^+} \Tt_{S_2}^+(X_2),
S_1^-\# T_2^+\right)
= \left(X_1\thsum{T_1 = S_2} X_2, S_1\# T_2\right).
$$
\end{lemma}
\proof{}  This is clear, since it is just a matter of cutting points out of
a neighborhood of $T_1$ in $X_1$  and moving them to a
neighborhood of $S_2^+$ in $\Tt_{S_2}^+(X_2)$.
Explicitly, using the definitions of thickening and thinning,
one can see that both sides are equal to
$$
\left(\Tt_{T_1}^-(X_1)\thsum{T_1^- = \Ga_k} W_{g,\eps}
\thsum{\Ga_{-k}=S_2} X_2,\, S_1^-\# F \# T_2\right)
$$
where $W_{g,\eps}$, $g=g_{T_1}$ is a ruled surface
with fiber $F$ and disjoint sections $\Ga_k,\Ga_{-k}$ where $k=-\io_{T_1}$.
\QED

The refined version of Lemma~\ref{le:resol} that we need is the following:

\begin{lemma}\label{le:resol2}
Given any triple $(X,S,T)$
and some  sufficiently small $\eps$, let $k=-\io_S$ and
consider the manifold $W_{g,2\eps}$, $g=g_S$, which has a section
$\Ga_k$ of area $A_S+(k+1)\eps$ where $A_S$ is the symplectic area of $S$.
Then
$$
\left(W_{g,2\eps} \thsum{\Ga_k= (S^-)^+} \Tt_{T^-}^{+}(\Tt_S^-(X)),\,
\Ga_{-k+2}\#(T^-)^+\right)
= \left(\Tt_{S+T}^{+}(X), (S+T)^+ \right)
$$
where the thickening and thinning are by the amount $\eps$.
Similarly, if $(g,k) = (g_T,\io_T+2)$ and the section $\Ga_{-k+2}$ has area
$A_T + (3-k)\eps$ then
$$
\left(\Tt_{S^-}^{+}(\Tt_T^- (X))\thsum{(T^-)^+ = \Ga_{-k+2}} W_{g,2\eps},\,
(S^-)^+\#\Ga_k\right)
 = \left(\Tt_{S+T}^{+}(X), (S+T)^+\right).
$$
\end{lemma}
\proof{}  As before, the second statement follows from the first by replacing
$k$ by $-k+2$ and interchanging $S,T$.  We begin by proving the
symplectic equivalence of the manifolds in the first statement, and
then show how the submanifolds are also identified.

Observe that the definitions immediately imply
$$
( W_{g,\eps}\thsum{ \Ga_k=S^-}\Tt_S^-(X) ,\, \Ga_{-k+2}\# T^-)  = (X, S+T).
$$
We need something more subtle.
Define
\begin{eqnarray*}
(X_3,S_3,T_3) & = & (\Tt_S^-(X),S^-,T^-)\\
(X_4,S_4,T_4) & = & (W_{g'\!,\eps}, \Ga_{k'}, F')
\end{eqnarray*}
where $g' = g_T, k' = -\io_T$ so that
$$
(\Tt_{T^-}^+(\Tt_S^-(X)),(S^-)^+)
 = (X_3\thsum{T_3=S_4} X_4,\, S_3\#T_4).
$$
The lemma then follows from the 4-fold sum rule (Proposition~\ref{prop:4sum})
once we express $W_{g,2\eps}$
appropriately as a sum of two ruled surfaces $X_1\# X_2$ each with fibers
$F^i$ of area $\eps$.  Namely, apply Lemma~\ref{le:eps0} using a
Hamiltonian function that has fixed point sets $\Ga_{-k+2},\Ga_{k-2}$
to define
\begin{eqnarray*}
(X_1, S_1,T_1) & = &  (W^1_{g,\eps}, F^1, \Ga_{k-2}^1)\\
(X_2, S_2,T_2) & = &  (W^2_{g,\eps}, \Ga_{-k+2}^2, \Ga_k^2)
\end{eqnarray*}
where the areas of $T_1 = \Ga_{k-2}^1$ and $S_2 = \Ga_{-k+2}^2$ both equal
$A_S+\eps$.
Then $S_1\# T_2$ has area  $A_S + (k+1)\eps$ and
$$
(X_1\,\#\,X_2)\,\#\,(X_3\,\#\,X_4) \,= \,
W_{g,2\eps} \thsum{\Ga_k= (S^-)^+} \Tt_{T^-}^{+}(\Tt_S^-(X)).
$$
On the other hand
$$
X_2\,\#\,X_3 \,=\, W^2_{g,\eps}\thsum{ \Ga_k=S^-}\Tt_S^-(X) \,=\, X
$$
and $X_4\#X_1$ is a ruled surface with fibers of area $\eps$
over a Riemann surface of genus
$g=g_S+g_T$.
The surface $S_2\#T_3 = \Ga_{-k+2}^2\#T^-$ has self-intersection
$\io_S+\io_T+2$ and is in the class of $[S]+[T]$, so we can choose it as the
representative $S+T$.
It is not hard to check that summing with $X_4\# X_1$ simply thickens
$X=X_2\#X_3$ along $S+T$ as desired, so
$$
(X_4\,\#\,X_1)\,\#\,(X_2\,\#\,X_3)\, =\, X_{S+T}^+.
$$

To see that the symplectomorphism of Proposition~\ref{prop:4sum}
identifies $\Ga_{-k+2}\#(T^-)^+$ and
$(S+T)^+$, notice that $(T^-)^+$ is in fact a section of $W_4$ disjoint
from $S_4$, and $\Ga_{-k+2}$ is a section of $W_1$ disjoint from $T_1$.
Therefore, when we sum along the fibers of $X_4,X_1$ we get another
ruled surface containing the connected sum of these sections, which is
a section disjoint from $S_4\#T_1$ and therefore corresponds to $(S+T)^+$.\QED

We are now ready to
prove that the associativity rule holds up to
symplectomorphism.
Roughly, the strategy is to thin the manifolds $X_i$ and
use the removed neighborhoods to create the necessary $X_4 = W_{g,2\eps}$ so
as to apply the 4-fold sum rule.
To do this correctly, we need to do some thickening as well.
Note that the order in which one thickens and thins does not matter, i.e.
$$
\Tt_{T^-}^+(\Tt_S^-(X)) = \Tt_{S^+}^-(\Tt_T^+(X)).
$$
\MS

\NI
{\bf Proof of Proposition~\ref{prop:assoc}}
We must show that
$$
(X_1\thsum{T_1 = S_2}X_2)\thsum{S_1\#T_2 = S_3 + T_3} X_3\;\; =\;\;
X_1\thsum{S_1+T_1=S_2\#T_3}(X_2\thsum{T_2=S_3}X_3)
$$
under the given hypotheses.
According to the above strategy,
consider the triples $(X'_i,S'_i,T'_i)$, $i=1,\cdots 3$ where we have
thickened and thinned by a sufficiently small $\eps$ to obtain
\begin{eqnarray*}
X'_1 & = & \Tt_{T_1^-}^+(\Tt_{S_1}^-(X_1)) \\
X'_2 & = & \Tt_{T_2^-}^-(\Tt_{S_2}^-(X_2)) \\
X'_3 & = & \Tt_{S_3^-}^+(\Tt_{T_3}^-(X_3))
\end{eqnarray*}
and surfaces $S'_i,T'_i$ which are the deformed $S_i,T_i$.
Then choose
$$
(X'_4,S'_4,T'_4) = (W_{g,2\eps},\Ga_{-k+2}, \Ga_k)
$$
where $k=-\io_S$ and the areas of $S'_4,T'_4$ equal
those of $T'_3,S'_1$.  This choice is possible since in both cases the
difference in the
areas of the two submanifolds is $2(1-k)\eps$.
The triples are admissible for the 4-fold sum, so the following
calculations prove the proposition.  We suppress the subscripts indicating
the submanifolds along which the sums are being performed when there
is no ambiguity.  Invoking Lemmas~\ref{le:eps} and~\ref{le:resol2}, we have
\begin{eqnarray*}
& & ( X'_1 \thsum{T'_1 = S'_2} X'_2 )
\thsum{ S'_1\#T'_2=S'_3\#T'_4}
( X'_3 \thsum{T'_3 = S'_4} X'_4 ) \\
& & \qquad \quad =
\left(\Tt_{T_1^-}^+(\Tt_{S_1}^-(X_1)) \,\#\,
\Tt_{S_2^-}^-(\Tt_{T_2}^-(X_2))\right)
\,\#\,
\left(\Tt_{S_3^-}^+(\Tt_{T_3}^-(X_3)) \,\#\, W_{g,2\eps}\right)\\
& & \qquad \quad =
\left(\Tt_{S_1}^-(X_1) \thsum{T^-_1=S^-_2}
\Tt_{T_2}^-(X_2)\right)
\thsum{S^-_1\#T^-_2 = (S_3+T_3)^+}
\Tt_{S_3+T_3}^+(X_3) \\
& & \qquad \quad =
(X_1\thsum{T_1 = S_2}X_2)\thsum{S_1\#T_2 = S_3 + T_3} X_3
\end{eqnarray*}
and
\begin{eqnarray*}
& & ( X'_4 \thsum{T'_4 = S'_1} X'_1 )
\thsum{ S'_4\#T'_1=S'_2\#T'_3}
( X'_2 \thsum{T'_2 = S'_3} X'_3 ) \\
& & \qquad \quad =
\left(W_{g,2\eps}\,\#\,
%_{\Ga_k=S'_1}
\Tt_{T_1^-}^+(\Tt_{S_1}^-(X_1))\right)
\,\#\,
%_{\{\Ga_{-k+2}\#T'_1 = S'_2\#T'_3\}}
\left(\Tt_{T_2^-}^-(\Tt_{S_2}^-(X_2)) \,\#\,
%_{T'_2=S'_3}
\Tt_{S_3^-}^+(\Tt_{T_3}^-(X_3))\right) \\
& & \qquad \quad =
\Tt_{S_1+T_1}^+(X_1) \thsum{(S_1+T_1)^+ = S_2^-\#T_3^-}
\left(\Tt_{S_2}^-(X_2) \thsum{T_2^-=S_3^-} \Tt_{T_3}^-(X_3)\right) \\
& & \qquad \qquad =
X_1\thsum{S_1+T_1=S_2\#T_3}(X_2\thsum{T_2=S_3}X_3)
\end{eqnarray*}.
\QED

\end{document}